\documentclass{elsarticle}
\usepackage{mathtools}
\usepackage{enumitem}
\usepackage{xcolor}
\usepackage{amsmath}
\usepackage{amsthm}
\usepackage{upgreek}
\usepackage{graphics}
\usepackage{makeidx}
\usepackage{color}
\usepackage{graphicx}
\usepackage{amssymb}
\usepackage{float}
\usepackage[english]{babel}
\usepackage{caption}
\usepackage{subcaption}
\usepackage{longtable}
\usepackage[utf8]{inputenc}
\usepackage[labelfont=bf]{caption}

\usepackage{amsthm}
\DeclareMathOperator{\Tr}{Tr}
\theoremstyle{definition}

\makeindex
\makeglossary
\begin{document}
\begin{frontmatter}
	\title{Modelling formation of stationary  periodic patterns in growing population of motile bacteria}
	
	\author{Valentina Bucur and Bakhtier Vasiev\corref{email}}
	
	\cortext[email]{bnvasiev@liverpool.ac.uk}
	
	\address{Department of Mathematical Sciences, Mathematical Sciences Building, University of Liverpool, Liverpool L69 7ZL, UK}
\begin{abstract}
Biological pattern formation is one of the most intriguing phenomena in nature. Simplest examples of such patterns are represented by travelling waves and stationary periodic patterns which occur during various biological processes including morphogenesis and population dynamics. Formation of these patterns in populations of motile microorganisms such as \textit{Dictiosteium dictiostelium} and \textit{E. coli} have been shown in a number of experimental studies.  Conditions for formation of various types of patterns are commonly addressed in mathematical studies of dynamical systems containing diffusive and advection terms. In this work, we do mathematical study of spatio-temporal patterns forming in growing population of chemotactically active bacteria. In particular, we perform linear analysis to find conditions for formation of stationary periodic patterns, and nonlinear (Fourier) analysis to find characteristics, such as amplitude and wavelength, of these patterns. We verify our analytical results by means of numerical simulations. 
\end{abstract}
\end{frontmatter}

\section{Introduction}
Pattern formation is one of the outcomes of biological self-organization and manifested by well known patterns, such as animal markings, and processes, such as embryogenesis. Growing bacterial communities form patterns which, for example, have an impact on biofilm formation \cite{James}. Pattern formation is a phenomenon which is widely observed in nature. Examples of non-biological pattern formation include formation of dunes in deserts, temperature hot spots in  autocatalytic chemical systems \cite{Gupta} and many other processes in various physical \cite{Smolin}, chemical \cite{Zaikin} and social \cite{vasiev2022} systems.

Turing, in his fundamental work \cite{Turing}, 
has pointed to two most representative examples of patterns, namely travelling waves and stationary periodic patterns. Formation of travelling waves in biological populations has been studied for many decades by means of Kolmogorov-Fisher equation \cite{Fisher,   Kolmogorov} and various models based on this equation \cite{Ablowitz, Fife}. One of the main futures of such models is an involvement of advection term which takes into account the directional motion of biological species in response to chemical signals (chemotaxis) \cite{Bramburger, Franz, Hosono, Keller1, Wang}.      

The most important outcome of Turing's work, presented in \cite{Turing}, was finding mechanism of stationary periodic pattern formation in chemical systems. Formation of this pattern was explained by so-called diffusion driven instability which nowadays is commonly referred as "Turing instability". Turing has looked at a system of differential equations describing the dynamics of two interacting substances, or morphogens. Commonly, in a well-mixed system (or in the absence of diffusion) the system has a stable homogeneous solution.   However, for certain types of interactions between the morphogens, i.e. in so-called "activator-inhibitor" system, homogeneous state becomes unstable when the diffusion of one of morphogens (inhibitor) is larger as compared with the diffusion of the other (activator). In this case, formation of non-homogeneous patterns such as stripes or spots takes place \cite{Budrene, Ks, Murray2}. This mechanism of pattern formation in the case of activator-inhibitor systems is commonly phrased as 'local self-enhancement and long-range inhibition' \cite{Meinhardt, Yamaguchi}.

Stationary periodic patterns are often observed in populations of chemotactically active microorganisms, which aggregate and form regularly spaced clusters. Chemotactic activity is expressed as motion along or against gradient of chemicals, produced externally or by microorganisms themselves.  Many theoretical studies have addressed the mechanisms of pattern formation in growing populations of such microorganisms represented by bacteria or amoeba \cite{Ks,Painter}. The general form of such models is given as 
\begin{equation} \label{general_model}
\begin{cases}
 \displaystyle \frac{\partial n}{\partial t}=D_n \nabla^2 n-\nabla \left( \chi(n,c) \nabla c \right) +f(n, c),\\[10pt]
  \displaystyle \frac{\partial c}{\partial t}=D_c \nabla^2 c + g(n,c),\\
    \end{cases}
\end{equation}
where the variable $n$ represents the density of microorganisms, and the variable $c$ - the concentration of chemotactic agent. Parameter $D_n$ is a diffusion coefficient describing the random motion of cells, $D_c$ - the diffusion constant of chemotactic agent, function $\chi(n,c)$ describes the chemotactic response of the microorganisms, and functions $f(n,c)$ and $g(n,c)$ describe the microbial and chemical kinetics, respectively. If $\chi(n,c) =0$ (no chemotaxis) the system \eqref{general_model} represents a reaction-diffusion model of the type of classical Lotka-Volterra model \cite{Murray1, Murray2}. If $\chi(n,c) \neq 0$ then we have an advection term in the first equation and the system \eqref{general_model} is a  reaction-diffusion-advection model. Furthermore, if $f(n,c)=0$, the system transfers into the well-known Keller-Segel model \cite{Ks}. 

Keller-Segel model was presented in 1970 and designed to describe the aggregation on \textit{Dictyosteium discoideum} amoeba mediated by chemotactic agent (cAMP, which at that time was known as acrasin) produced by the amoeba themselves \cite{Ks}. The original model consisted of four equation, but using quasi-steady-state assumption (i.e. by replacing 'fast' differential equations by algebraic relationships) it was reduced to a system of two reaction-diffusion-advection equations. The analysis was performed for the simplest version of the model where it was assumed that the chemotactic sensitivity is constant ($\chi(n,c)=\chi_0$), amoeba do not die and reproduce, ($f(n,c)=0$), and the chemical agent is produced by the amoeba and degrade at a constant rate, $g(n,c)=hn-pc$, where $h$ is the production rate per amoebae and $p$ is the degradation rate. Under these assumptions, the model \eqref{general_model} is given by the following system:  
\begin{equation*}
\begin{cases}
 \displaystyle   \frac{\partial n}{\partial t} = D \nabla^2 n-\chi_0 \nabla^2 c, \\[5pt]
 \displaystyle   \frac{\partial c}{\partial t}= \nabla^2c +hn-pc,\\
\end{cases}
\end{equation*}
where $D=D_n/D_c$ is the ratio of diffusion constants of the amoeba and chemical and the spatial variable is rescaled in a way that the diffusion of the chemical $c$ is set to one. In the absence of chemotaxis ($\chi=0$) the system evolves towards the stable steady state $(n,c)=(n_0,c_0$) defined by the initial amount of amoeba, so that $n_0$ is an average initial density of amoeba and $c_0=n_0*h/p$. It was shown in \cite{Ks} that in the presence of chemotaxis the homogeneous state ($n_0,c_0$) becomes unstable if
\begin{equation*}
\displaystyle \frac{\chi_0 p}{Dh}+\frac{n_0h}{p}>1.
\end{equation*}
This condition implies that an increase in the chemotactic activity $\chi$ breaks down the stability of homogeneous state and the amoeba aggregate into clusters forming stationary periodic patterns. An assumption that the chemotactic response doesn't depend on the density of cells and on the concentration of chemotactic agent (that is, $\chi(n,c)=\chi_0$) is oversimplification which was made in \cite{Ks} to make the mathematical analysis easier. It is reasonable to assume that the chemotactic response is proportional to the density of cells, and there is experimental evidence indicating that it is reduced with an increase of the concentration of chemotactic agent. Our analysis of the modification of model, where the chemotactic term was assumed having the form  $\displaystyle \chi(n,c)=\chi_0\frac{n}{c}$, have shown that the aggregation takes place if \mbox{ } $\displaystyle  \frac{\chi_0}{D}+\frac{n_0h}{p}>1$. This result doesn't alter the main conclusion drawn from the previous version of the model that the chemotactic response should be large enough to allow formation of stationary periodic patterns. Other modifications of Keller-Segel model (involving different forms of chemotactic term and/or bacterial and chemical kinetics) have been considered by a number of researchers and some on them are listed in \cite{Painter}. 

\textit{E. coli} bacteria is another biological species which exhibits chemotactic activity.  Patterns formed in \textit{E. coli} colonies have been studied experimentally in great details \cite{Budrene1,Budrene}. It was shown that the waves of expansion of growing \textit{E. coli} colonies transform into stationary spots formed by clusters of dead/quiescent bacteria. It was also shown that these patterns form because \textit{E. coli} respond to chemoattractant which they secrete themselves. It was established that the observed patterns of cell aggregation depends on nutrient consumption, cell proliferation, excretion of attractant and chemotactic motility \cite{Budrene1}. A number of models have been developed to describe patterns formed in growing \textit{E. coli} colonies \cite{Mimura,Polezhaev}. 
The earliest models reproducing the formation of spots in  \textit{E. coli} colonies were presented in \cite{Tyson}. One of these models (so-called 'liquid model' which is designed to reproduce patterns formed by bacteria in liquid environment) is very similar to Keller-Segel model and given by the following two equations:
\begin{equation*}
\begin{cases}
\displaystyle \frac{\partial n}{\partial t}=D \nabla^2n-\chi_0\nabla\left(\frac{n}{(1+c)^2}\nabla c\right), \\[7pt]
\displaystyle \frac{\partial c}{\partial t}=\nabla^2c+\omega \frac{n^2}{\nu+n^2}.
\end{cases}
\end{equation*}
In this model the death and proliferation of cells is again neglected ($f(n,c)=0$). The chemotactic sensitivity is given as $\displaystyle \chi(n,c)=\chi_0\frac{n}{1+c^2}$, that is, it is assumed to be proportional to the density of cells and the dependence on concentration of the chemotactic agent is modified in a way that the absence of chemical doesn't imply the division by zero. The decay of the chemical is neglected while its production is given by the Hill function $\displaystyle g(n,c)=\omega \frac{n^2}{\nu+n^2}$. Analysis of this model performed in \cite{Tyson} has shown that for sufficiently large $\chi_0$ the homogeneous state becomes unstable with formation of periodic patterns for which, however, the spacial periodicity changes over on time. Thus, the patterns obtained using models presented in \cite{Tyson} are not stationary. Formation of stationary periodic patterns was shown later by other researchers in a framework of models for dynamics of \textit{E. coli} population represented by larger systems of differential equations \cite{Mimura,Polezhaev}. 

While mechanisms and conditions of periodic pattern formation have been addressed in many studies \cite{Turing, Ks, Murray1} their properties, such as wavelength and amplitude, where studied in a much lesser extend. The estimation of wavelength of patterns is commonly performed on the basis of linear instability analysis and considered to be given by the wave number of most unstable mode (see Chapter 2 of \cite{Murray2}). Most recently properties of Turing patterns formed in a class of two-variable reaction-diffusion models, with kinetics terms given by polynomial functions, were addressed in \cite{Buceta}. Amplitude of the patterns were estimated in assumption that they have a sinusoidal form, that is, that the solution of \eqref{general_model} (in one-dimensional case) is given as:
\begin{equation*}
\begin{cases}
n(x,t)=n_0+\alpha \cos(kx), \\
c(x,t)=c_0+\beta \cos(kx),  
\end{cases}
\end{equation*} 
where ($n_0$, $c_0$) represents homogeneous steady state of the system, $k$ is the wave number of most unstable mode found from linear stability analysis and $\alpha$ and $\beta$ give amplitude of Turing pattern emerging in the model. The estimate of the amplitude, obtained this way, has matched reasonably well numerical results presented in this research. However, it is evident that more accurate estimates of the wavelength and amplitude of Turing patterns should be based on nonlinear analysis of these patterns.  

In this paper we study properties of periodic patterns forming in a representative version of model \eqref{general_model} which is relatively simple, justifiable biologically, and allows such pattern formation. The choice of the model was made on the basis of analysis of patterns forming in a range of models of type   \eqref{general_model} with non-zero advection term. We use Fourier series for analysis of patterns obtained in numerical simulations as well as for analytical estimation of the amplitude and wavelength of Turing patterns emerging in the model. For greater accuracy we perform our analytical studies for a small domain which can  contain only half (or one) spike. We also study how the amplitude and wavelength of Turing patterns depend on the values of model parameters.


\section{Model}

In this section, we introduce the model of growing population of motile bacteria which we will analyse throughout this paper. This model is represented by a system of two partial-differential equations which include diffusion, reaction and advection terms. It is known that such systems allow formation of travelling wavefronts \cite{Franz, Li} as well as Turing patterns \cite{Painter, Murray1, Murray2}. 

We look for Turing patterns forming in one-dimensional domain (of size $L$) under no-flux boundary conditions. We have performed preliminary numerical simulations using various models of type \eqref{general_model} and based on the results of these simulations we have chosen the following system: 
\begin{equation*}
\begin{cases}
\displaystyle \frac{\partial n}{\partial t}=D_n \frac{\partial ^2 n}{\partial x^2}-\tilde{\chi}\frac{\partial}{\partial x}\left(n\frac{\partial c}{\partial x}\right)+r_0n\left(1-\frac{n}{k}\right),\\[10pt]
\displaystyle \frac{\partial c}{\partial t}= D_c\frac{\partial ^2 c}{\partial x^2} +hn-pc,\\
\end{cases}
\end{equation*}
This system is a version of the system \eqref{general_model} and describes the dynamics of population of bacteria, whose density is given by $n$, which respond chemotactically to a chemical agent of concentration $c$, produced by bacteria themselves. Here the chemotactic sensitivity of cells is assumed to be proportional to $n$ and takes the form $\chi(n,c)=\tilde{\chi} n$. Also, cells reproduce according to logistic law with growth rate $r_0$ and carrying capacity $k$: $f(n,c)=r(1-r/k)$. The kinetics of chemical is assumed to be linear , $g(n,c)=hn-pc$ (like in Keller-Segel model \cite{Ks}), that is, it is being produced by cells with rate $h$ and decays with rate $p$. This system can be nondimensionalised using the substitutions $\tilde{t}=tp$, $\displaystyle \tilde{x}=x\sqrt{\frac{p}{D_c}}$, $\displaystyle \tilde{n}=\frac{n}{k}$, $\displaystyle \tilde{c}=\frac{cp}{hk}$. This gives (after dropping the tildes): 
\begin{equation} \label{nondimensional}
\begin{cases}
\displaystyle \frac{\partial n}{\partial t}=D\frac{\partial ^2 n}{\partial x^2}-\chi_0 \frac{\partial}{\partial x}\left(n\frac{\partial c}{\partial x}\right)+rn(1-n),\\[10pt]
\displaystyle \frac{\partial c}{\partial t}= \frac{\partial ^2 c}{\partial x^2} +n-c,\\
\end{cases}
\end{equation}
with new constants $\displaystyle D=\frac{D_n}{D_c}$, $\displaystyle \chi_0=\frac{\tilde{\chi}hk}{D_cp}$ and $\displaystyle r=\frac{r_0k}{p}$ While the parameters $D$ and $r$ are essentially non-negative, $\chi$ can be positive (chemoattraction) or negative (chemorepulsion).

In the absence of diffusion and chemotaxis, the system is said to be well-mixed, meaning that both variables, $n$ and $c$, don't depend on spatial variable $x$. In this case, the first equation is detached from the second and describes the logistic growth of population of bacteria. The system has two steady states $(n, c)=(0,0)$ and $(1,1)$, which are unstable and stable respectively \cite{Murray1}. After putting back the diffusion term ($D \ne 0$, $\chi_0=0$) the first equation transfers into Kolmogorov-Fisher equation which is known for travelling front solutions \cite{Kolmogorov, Fisher}. The system now allows formation of travelling wavefronts between the steady states $(n, c)=(0,0)$ and $(n,c)=(1,1)$, which are now stable and unstable, respectively \cite{Murray1}. 

\begin{figure}[h]
\begin{subfigure}{.32\textwidth}
\includegraphics[scale=0.22]{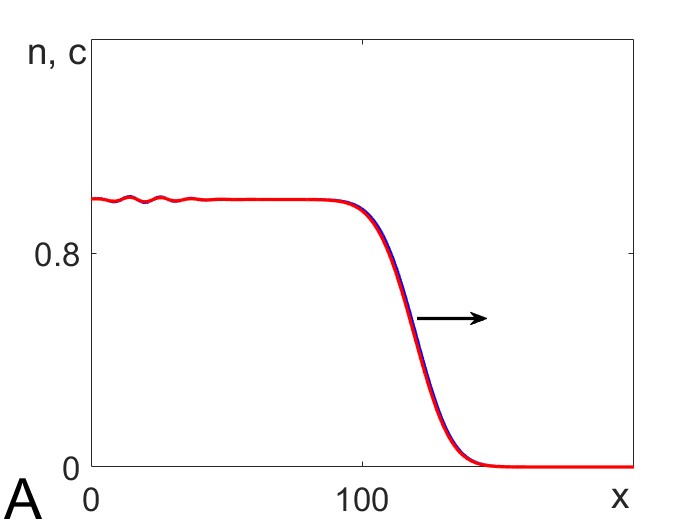}
\end{subfigure}
\begin{subfigure}{.32\textwidth}
\includegraphics[scale=0.22]{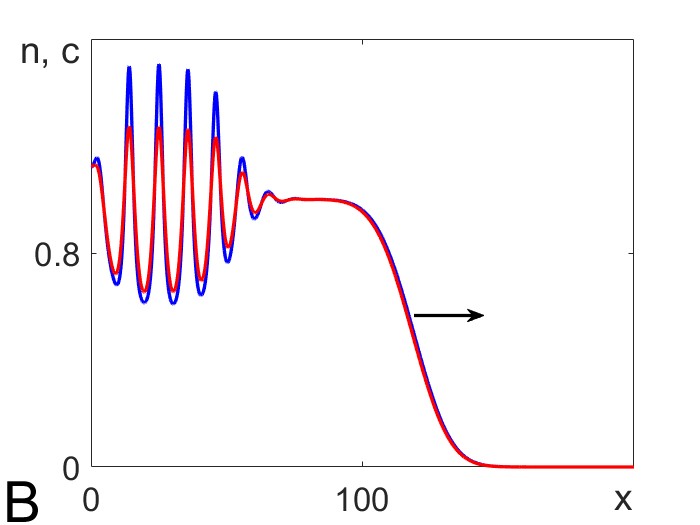}
\end{subfigure}
\begin{subfigure}{.32\textwidth}
\includegraphics[scale=0.22]{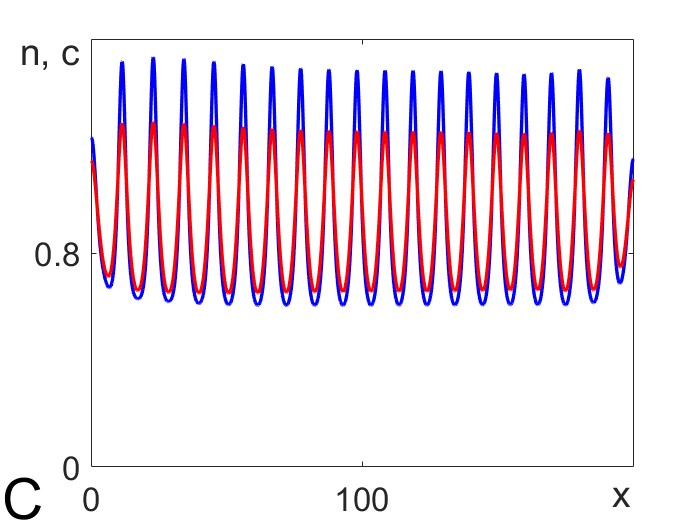}
\end{subfigure}
\caption{\em{\textbf{Numerical simulation of travelling front and Turing patterns forming in the system (\ref{nondimensional})}. Profiles of cell density, $n$, (blue line) and concentration of chemotactic agent, $c$, (red line) are shown. Solutions are obtained in the domain $0 \le x \le 200$, $t \ge 0$.  Initial conditions: $n(x,0)=0$ if $x>2$ and $n(x,0)=1$ if $0 \le x \le 2$; $c(x)=0$  $\forall x$. \textbf{A}: formation of travelling front in a system when the chemotaxis is weak ($\chi_0=1.7$). \textbf{B, C}: formation of Turing pattern in a system when the chemotaxis is strong ($\chi_0=1.9$). Other model parameters: $D=1$, $r=0.1$. Profiles are shown at time  $T=200$  in panels \textbf{A} and \textbf{B}, and $T=800$  in panel \textbf{C}.}}
\label{fig_TP_TW}
\end{figure}

When the system includes the flow of cells due to chemotaxis ($\chi_0 \ne 0$) then the two equations in \eqref{nondimensional} are coupled and the known solutions of this system include stationary periodic, or Turing, patterns \cite{Murray1}. When the chemotaxis is weak ($\chi_0$ is small), we have a model, which, similarly to the Fisher-Kolmogorov equation, allows only formation of travelling wavefront (see Fig.\ref{fig_TP_TW}A). However, if the chemotactic flow is strong enough, that is, the value of $\chi_0$ is above certain value, then the formation of stationary periodic patterns given by periodic functions, $n(x)$ and $c(x)$, can be observed (see Fig.\ref{fig_TP_TW}B,C). 

\section{Conditions for formation of Turing patterns}

Our numerical simulations, illustrated in Fig.\ref{fig_TP_TW}, have indicated that Turing patterns can appear only if the chemotactic sensitivity, $\chi_0$, is above some threshold value. This threshold value can be found analytically by applying the standard linear analysis technique, which is well known from the literature \cite{Ks, Murray2}. It is known that Turing patterns form due to so-called Turing instability when a steady state, which is stable in the well-mixed system, becomes unstable in the full reaction-diffusion-advection system. This means that the equilibrium point $(n_0, c_0)=(1, 1)$ in system \eqref{nondimensional}, which is stable if $D=\chi_0=0$, should become unstable for some non-zero values of $D$ and $\chi_0$. 

Investigation of stability of steady state $(n_0,c_0)$ involves linearisation of the system \eqref{nondimensional} at this state and checking the dynamics of perturbations, which are given in the form $e^{\lambda t}\cos(kx)$. Parameter $k$ represents a perturbation mode and can have any value from the set, given as $k=i\pi/L$, with $i$ taking any integer value so that the perturbation satisfies no-flux boundary conditions. It appears that these perturbations represent solutions of the linearised system provided that $\lambda$ is an eigenvalue of the characteristic matrix, which is, for the system \eqref{nondimensional} at the steady state $(1,1)$, given as
\begin{equation} \label{detM}
M=
\begin{pmatrix}
-Dk^2 -r & \chi_0 k^2\\
1  & -k^2-1
\end{pmatrix}.
\end{equation}
The steady state $(1,1)$ is stable if both eigenvalues of this matrix are negative, which is known to require $\det M>0$ and $\Tr M<0$ \cite{Murray2}. We can see that $\Tr M=-Dk^2 -r-k^2-1$ is always negative, so for the system to be driven unstable by perturbation, we require $\det M<0$, which leads to the following condition for Turing instability,
\begin{equation*} 
\det M=Dk^4+k^2(D +r - \chi_0 )+r<0.
\end{equation*}
From this, a necessary but not sufficient condition is 
\begin{equation} \label{cond2}
D +r-\chi_0 <0
\end{equation}
since if $D +r - \chi_0>0$ then $\det M>0$ for any $k$.The second necessary condition for initiation of Turing patterns comes from the fact that the equation $\det M=0$ should have real roots, so that $\det M<0$ between these two roots. Clearly, the roots are given as
\begin{equation} \label{wavenumber}
\displaystyle k^2_{1,2}= \displaystyle \frac{-(D+r-\chi_0) \pm \sqrt{(D+r-\chi_0)^2-4Dr}}{2D},
\end{equation}
and therefore we require
\begin{equation*}
(D+r-\chi_0)^2-4Dr>0.
\end{equation*}
Taking into account \eqref{cond2} we get
\[\chi_0  -D-r>2\sqrt{Dr}. \]
After rearranging terms in this inequality and, by introducing new notation $R_T$, we can state the condition for Turing instability in the form: 
\begin{equation} \label{m3_rt}
R_T=\frac{\chi _0}{D  +r +2\sqrt{D r}}>1.
\end{equation}

Looking again at Figure \ref{fig_TP_TW}, we see that in panel A, the factor of instability is $R_T=0.98<1$, while in panels B and C, $R_T=1.1>1$. Thus, results shown in Fig.\ref{fig_TP_TW} confirm that if $R_T<1$ there observed travelling waves rather than Turing pattern, while if $R_T>1$ formation of Turing patterns is observed. 

Using the above technique we have performed stability analysis for a number of models, which are represented by system \eqref{general_model} but differ by functions $\chi(n,c)$ and $f(n,c)$. Following \cite{Ks} we assumed that the chemical kinetics is linear and didn't alter it across models, i.e. $g(n,c)=n-c$ in all models. List of the analysed models together with the obtained formulas for instability factor, $R_T$, is shown in Table \ref{tab_models}. Note that travelling wave solutions form if $R_T<1$ (and $f(n,c) \neq 0$), while if $R_T>1$ then formation of Turing patterns should be observed.

\begin{figure}[H]
\begin{subfigure}{.32\textwidth}
\includegraphics[scale=0.22]{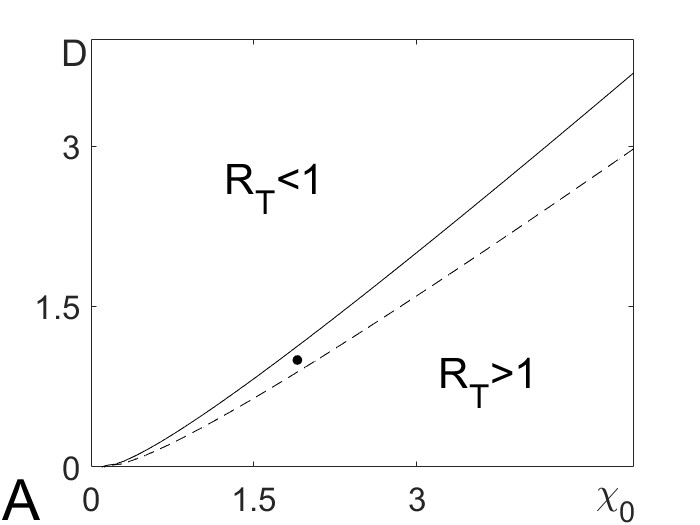}
\end{subfigure}
\begin{subfigure}{.32\textwidth}
\includegraphics[scale=0.22]{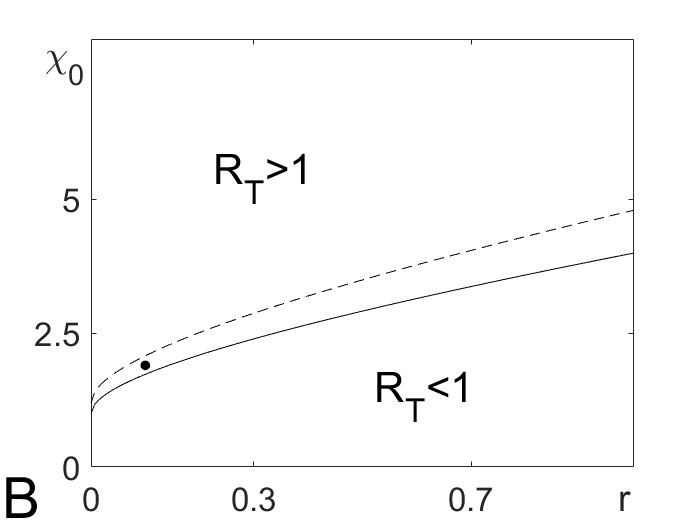}
\end{subfigure}
\begin{subfigure}{.32\textwidth}
\includegraphics[scale=0.22]{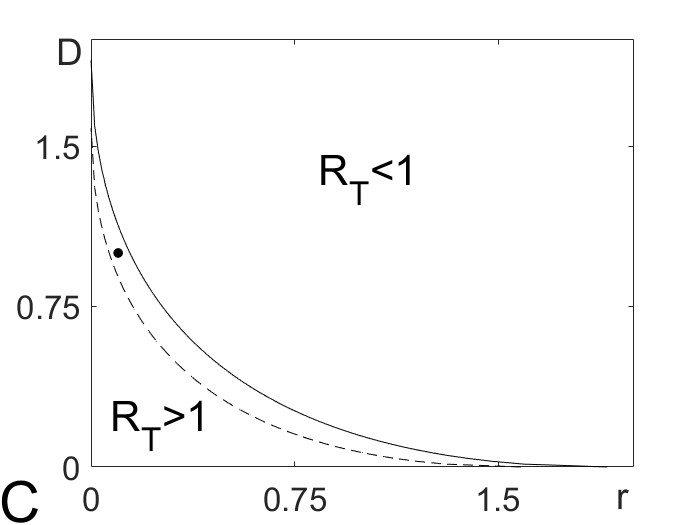}
\end{subfigure}
\caption{\em{\textbf{
Domains $R_T>1$ (where Turing instability takes place)  in parameter planes for  models $\mathbf{M3}$ (solid) and $\mathbf{M7}$ (dashed).} Plots $R_T=1$ are drown according to the formulas for $R_T$ given in Table \ref{tab_models} in an assumption that one of parameters has a fixed (default) value, which are $D=1$, $\chi_0=1.9$,  $r=0.1$. Dashed line is plotted for $\nu=0.2$ (model $M7$). The dot indicates position of the point with coordinates given by default values of parameters.}}
\label{ratio_implicit}
\end{figure}

\begin{table}[h]
\begin{center}
\begin{tabular}{| l | l | l | l |} 
    \hline
  & $\chi(n,c)$ & $f(n, c)$ & $R_T$ \\[7pt] \hline
$M_1$ & $\chi_0$ & 0 & $\displaystyle \frac{\chi_0}{D}$ \\[7pt] \hline
$M_2$ & $\chi_0n$ & 0 & $\displaystyle \frac{\chi_0n_0}{D}$ \\[7pt] \hline
$M_3$ & $\chi_0n$ & $rn(1-n)$ & $\displaystyle \frac{\chi_0}{D+r+2\sqrt{Dr}}$ \\[7pt] \hline
$M_4$ & $\displaystyle \chi_0\frac{1}{c}$ & 0 &  $\displaystyle \frac{\chi_0}{Dn_0}$\\[10pt] \hline
$M_5$ & $\displaystyle \chi_0\frac{1}{c}$ & $rn(1-n)$ & $\displaystyle \frac{\chi_0}{D+r+2\sqrt{Dr}}$ \\[10pt] \hline
$M_6$ & $\displaystyle \chi_0\frac{n}{c}$ & $rn(1-n)$ & $\displaystyle \frac{\chi_0}{D+r+2\sqrt{Dr}}$ \\[10pt] \hline
$M_7$ & $\displaystyle \chi_0\frac{n}{c+ \nu}$ & $rn(1-n)$ & $\displaystyle \frac{\chi_0}{(\nu + 1)(D+r+2\sqrt{Dr})}$ \\[10pt] \hline
$M_8$ & $\displaystyle \chi_0\frac{n}{(c+ \nu)^2}$ & $rn(1-n)$ &  $\displaystyle \frac{\chi_0}{(\nu + 1)^2(D+r+2\sqrt{Dr})}$ \\[10pt] \hline
$M_9$ & $\displaystyle \chi_0\frac{n+\nu_1}{c+ \nu_2}$ & $rn(1-n)$ &  $\displaystyle \frac{\chi_0(1+\nu_1)}{(\nu_2 + 1)(D+r+2\sqrt{Dr})}$ \\[10pt] \hline
    \end{tabular}
\end{center}
\caption{\em{\textbf{Results of Turing instability analysis for nine models represented by system \eqref{general_model}}. Model notations are given in the first column. Functions $\chi(n,c)$ and $f(n,c)$ corresponding to each model are given in columns 2 and 3. Calculated formula for instability factor, $R_T$, is given in column 4. In models M2 and M4, where $f(n,c)=0$, the instability factor, $R_T$, depends on $n_0$, which is an average value of $n$ in the domain and determined by the initial conditions.}}
\label{tab_models}
\end{table}
Most of the models listed in Table \ref{tab_models} have been introduced earlier (see \cite{Painter}), however, the conditions for Turing instability for most of these models have not been stated so far. These conditions are qualitatively similar for all considered models and this is illustrated by domains where these conditions satisfied (in parameter planes for models $M3$ and $M7$) in Fig. \ref{ratio_implicit}. We see from this figure that for Turing instability to take place the chemotactic sensitivity, $\chi_0$ should be strong enough, while the diffusion of chemotactic agent, $D$ and reproduction rate $r$ are small enough. 

\begin{figure}[H]
\begin{subfigure}{.32\textwidth}
\includegraphics[scale=0.22]{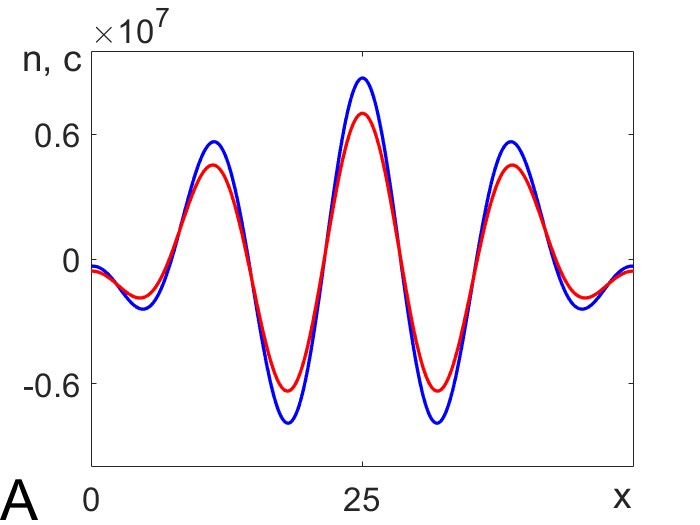}
\end{subfigure}
\begin{subfigure}{.32\textwidth}
\includegraphics[scale=0.22]{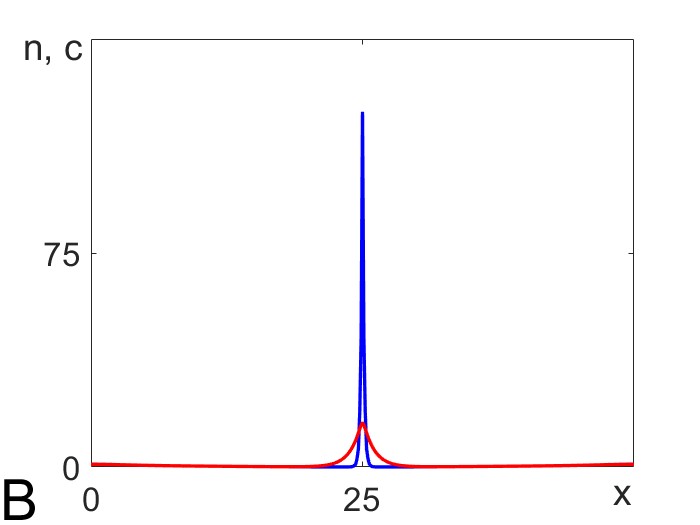}
\end{subfigure}
\begin{subfigure}{.32\textwidth}
\includegraphics[scale=0.22]{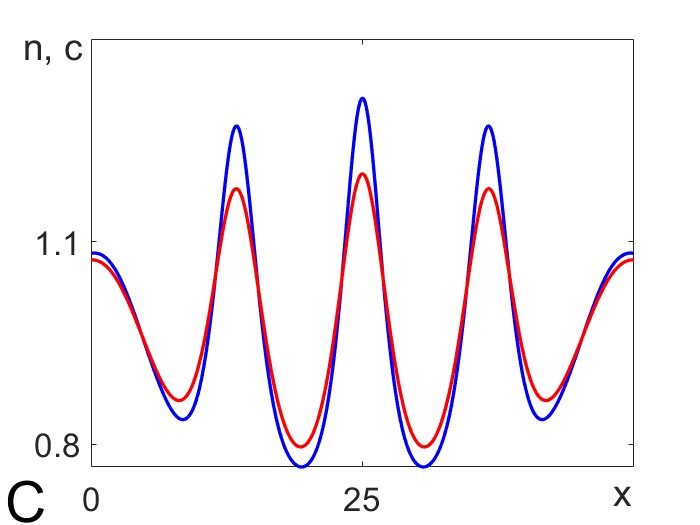}
\end{subfigure}
\caption{\em{\textbf{Numerical simulations of patterns forming in models $M_1$ (panel A), $M_2$ (panel B) and $M_3$ (panel C).} Profiles of cell density, $n$, (blue line) and concentration of chemotactic agent, $c$, (red line) at time $T=500$ in the domain of size $L=50$ shown in all three panels. Initial condition $n=c=0$ in (\textbf{A}) and (\textbf{C}), and $n=1$ and $c=0$ in (\textbf{B}). Patterns were initiated using stimulus applied to the centre of the domain: $n=1.1$ for $24<x<26$ at time $t=0$. Model parameters: $D=1$, $\chi_0=1.8$ and $r=0.1$ (in $M_3$).}} 
\label{three_models}
\end{figure}

We have reproduced numerically Turing patterns forming in all models given in Table \ref{tab_models} aiming to find the simplest model, which allows formation of Turing patterns, such that its shape is within physical constrains. Here is a summary of our  numerical results: 
\begin{itemize}
	\item $M_1$ is the simplest model where the chemotactic sensitivity is assumed to be constant, and no cell death/proliferation takes place \cite{Ks}. Numerical simulations show that Turing patterns in this model evolve towards infinite amplitude and involve negative values for density of cells and concentration of chemical (see Fig.\ref{three_models}A), which is physically impossible. So, we can conclude, that $\chi(n,c)$ has to depend on at least one of the two variables, $n$ and $c$.
	\item  In $M_2$ the chemotactic sensitivity is assumed to be proportional to the concentration of cells, $\chi(n,c)=\chi_0 n$, but there is still no cell death/proliferation are taking place. When the instability condition is met we observe only one sharp spike forming in response to stimulus (see Fig. \ref{three_models}B). Thus, no Turing pattern can be observed in this model. 
	\item  $M_3$ is $M_2$ with added logistic growth of cells, $f(n,c)=rn(1-n)$. In this model Turing patterns of finite amplitude and non-negative values of $n$ and $c$ are observed. Thus, this is a simplest model which allows reproduction of Turing pattern of physically justifiable shape. 
	\item In $M_4$ the chemotactic sensitivity is inversely proportional to the concentration of cells and there is no cell proliferation. Here numerical simulations end up with division by zero and abruptly stop after short period of time.
	\item $M_5$ is $M_4$ with added logistic growth of cells. Simulations can be ran for longer periods of time, however division by zero is still taking place. 
	\item In $M_6$ $\chi(n,c)$ is proportional to $n$ and inversely proportional to $c$. Eventually, division by zero is not taking place here and we obtain appropriately shaped Turing patterns.
	\item  The main feature of models $M_{7, 8, 9}$ is that $\chi(n,c)$ is proportional to linear function of $n$ and inverse proportional to linear function, or a square of linear function of $c$. The latter property lets us to safely avoid division by zero (see \cite{Painter, Tyson}). In all these cases, we get reasonably shaped Turing patterns, which can be sharper or smoother for particular models (i.e. $M_7$, $M_8$ or $M_9$). 
\end{itemize}

Our main conclusion after analysis and numerical simulations presented in this chapter that model $M_3$ is the simplest model where periodic stationary patterns, which do not violate any physical constrains, are observed. Model $M_3$ (which is given by the system \eqref{nondimensional}) will be used in the rest of this paper for analysis of properties of Turing patterns such as their amplitude and wavelength.

\section{Estimation of the wavelength of Turing patterns}

In the previous sections, we have found conditions for formation of Turing patterns (given, for example, by Eq. \eqref{cond2} for the model \eqref{nondimensional}) and have presented numerical simulations (see Fig.\ref{fig_TP_TW}) to illustrate the validity of this estimation. Our next task is to analyse properties of Turing patterns, which are given by their wavelength and amplitude. The amplitude is defined as the difference between maximal and minimal values of cell density, $n$, in the pattern,  while the wavelength - as a distance representing spatial periodicity of Turing pattern. As an alternative to the wavelength we will also consider a characteristic length which is a distance between the consecutive maxima and minima and twice smaller than the wavelength. 

In this section we consider the wavelength of Turing pattern and will estimate it on the basis of linear stability analysis presented in the previous section. Turing pattern is a periodic pattern, and as a such can be represented as a sum of cosinusoidal profiles, each given as $a(k)\cos(kx)$, where $k=i\pi/L$, $i$ is an integer number and $L$ is the domain size. Integer $i$ should be such that $k$ belongs to the interval $k \in (k_1, k_2)$ where $k_1$ and $k_2$ are given by Eq. \eqref{wavenumber}. Wavenumbers $k$ in this interval correspond to unstable modes, each having the wavelength $\Lambda=\pi/k$. It is reasonable to assume that Turing pattern is predominantly defined by the most unstable mode. Popular way to choose this mode is to take the one corresponding to the average of $k_1^2$ and $k_2^2$ given by Eq.\eqref{wavenumber} \cite{Murray2, Buceta}, which is: 
\begin{equation*}
\displaystyle k^2_{av}= \displaystyle \frac{\chi_0-D-r}{2D}.
\end{equation*}
The wavelength corresponding to this mode is given as:
\begin{equation} \label{lamb1}
\displaystyle \Lambda=\pi\sqrt{\frac{2 D}{\chi_0-D-r}}.
\end{equation}

More accurate way to estimate the most unstable mode is to find a mode for which the characteristic matrix \eqref{detM} has the maximal positive eigenvalue. In terms of trace and determinant, the two eigenvalues of the characteristic matrix, $M$, are 
$$\displaystyle \lambda_{1,2}=\frac{tr M \pm \sqrt{tr M ^2-4\det M}}{2},$$
which gives
\begin{equation*}
\displaystyle \lambda_{1,2}= \frac{-Dk^2 -k^2 -r-1}{2} \pm \frac{\sqrt{k^4(D^2-2D+1)+k^2(2D(r-1)+4\chi_0-2r+2)+1-2r}}{2}.
\end{equation*}
Clearly, $\lambda_1<0$ always, so we are interested in the value of $k^2$ which corresponds to the largest value of $\lambda_2$. This means that we are looking for the value of $k^2$, at which the derivative of $\lambda_2$ with respect to $k^2$ is zero: 
\begin{equation}\label{lamb2}
\displaystyle \frac{d \lambda_2}{d k^2}=\frac{-D-1}{2}+\frac{\tau(1-D)^2+(2\chi_0 +Dr-D-r+1)}{2\sqrt{\tau^2(1-D)^2+2\tau(2\chi_0 +Dr-D-r+1)-2r+1}}=0.
\end{equation}
Using Matlab code, we have found the value of $k^2$ (for fixed values of parameters $D$, $\chi_0$ and $r$) satisfying the above equation, and hence calculated the corresponding wavelength, $\Lambda=\pi/k$. 
\begin{figure}[H]
	\begin{subfigure}{.32\textwidth}
		\includegraphics[scale=0.22]{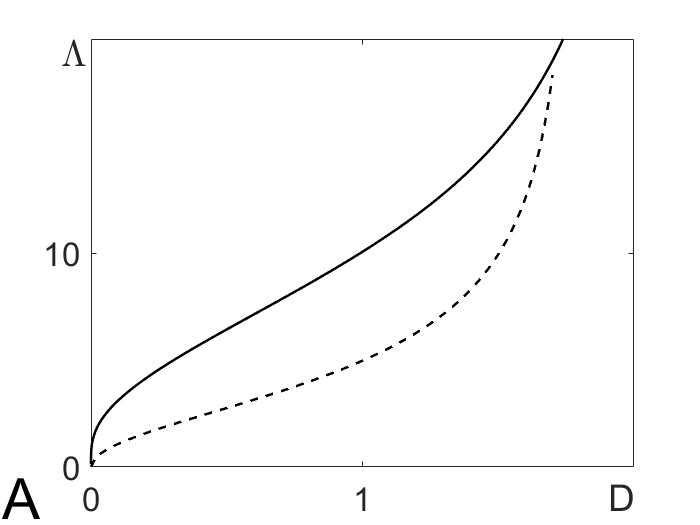}
	\end{subfigure}
	\begin{subfigure}{.32\textwidth}
		\includegraphics[scale=0.22]{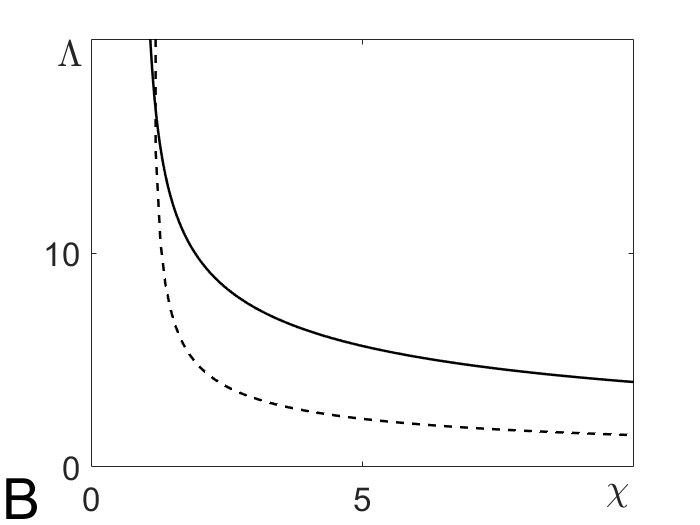}
	\end{subfigure}
	\begin{subfigure}{.32\textwidth}
		\includegraphics[scale=0.22]{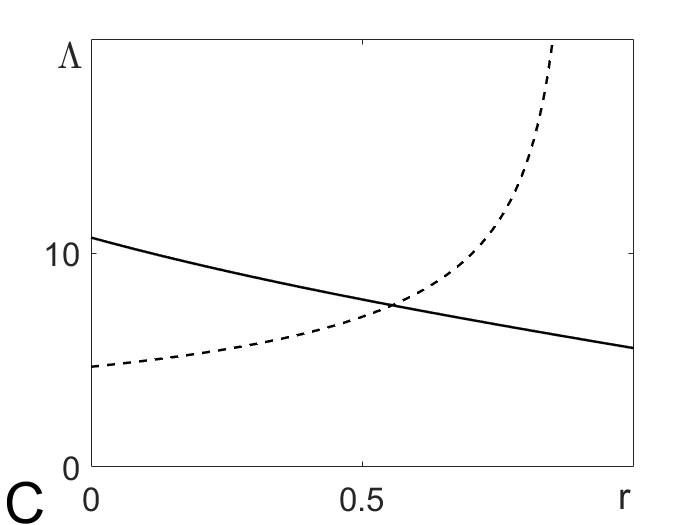}
	\end{subfigure}
	\caption{\em{\textbf{
Dependence of the wavelength of Turing pattern, $\Lambda$, on parameters of model \eqref{nondimensional} found according to the equation \eqref{lamb1} (dashed) and the implicit formula \eqref{lamb2} (solid).} Default values of model parameters: $D=1$, $\chi_0=1.9$ and $r=0.1$. 
	}}
	\label{linear_approx}
\end{figure}

To complete the presented linear analysis we have checked the dependence of the wavelength, $\Lambda$,  on the values of model parameters for the system \eqref{nondimensional} according to formula \eqref{lamb1} or implicit formula \eqref{lamb2}. The corresponding plots are shown in Fig. \ref{linear_approx}. One can see that according to both formulas the wavelength is monotonically increasing function of diffusion, $D$ (panel A), and decreasing function of chemotactic sensitivity, $\chi_0$ (panel B). Vertical asymptotes in panels A and B confirm that Turing patterns can't be observed, when the diffusion is too fast, or the chemotactic sensitivity is too weak (see the condition given by Eq.\eqref{m3_rt} and domains for Turing instability in Fig. \ref{ratio_implicit}). While we observe certain quantitative differences between plots in Panels A and B, plots on Panel C differ qualitatively: the wavelength of most unstable mode in Turing pattern is a decreasing function of the growth rate, $r$, according to formula \eqref{lamb2} and an increasing function according to the formula \eqref{lamb1}. Our numerical simulations (shown in Section 7) indicate that the dependence given by the formula \eqref{lamb2} is not only correct qualitatively but also quite accurate quantitatively. 

We note that Equations \eqref{lamb1} and \eqref{lamb2} only give approximations of the wavelength for Turing pattern based on the linear analysis of the system \eqref{nondimensional}. These approximations will later be compared with results obtained by means of nonlinear analysis and numerical simulations. 


\section{Fourier series for numerically simulated profiles}

Turing patterns are stationary periodic patterns and as such can naturally be represented as Fourier series. Pattern forming in a domain of size $L$ under no-flux boundary conditions is given as a superposition of cosines, that is in the case of the system \eqref{nondimensional} we have:
\begin{equation} \label{sol}
\begin{cases}
\displaystyle n(x)=\sum_{i=0}^M \alpha_i\cos\frac{i\pi x}{L},\\[10pt]
\displaystyle c(x)=\sum_{i=0}^M \beta_i\cos\frac{i\pi x}{L}.
\end{cases}
\end{equation}
The coefficients $\alpha_i$ and $\beta_i$ define the amplitudes of mode $i$ for the variables $n$ and $c$ (the amplitudes are given as $2\alpha_i$ and $2\beta_i$). For smooth profiles these coefficients quickly tend to zero as $i$ increases and therefore we can truncate the series by taking the first $M$ terms, with $M$ to be carefully defined. For known profile, $n(x)$, the coefficients are found using the formulas
\begin{equation*}
\alpha_0=\frac{1}{L}\int_{0}^{L}n(x)dx, 
\end{equation*}
and for $i>0$
\begin{equation} \label{alphas}
\alpha_i=\frac{2}{L}\int_{0}^{L}n(x)\cos\frac{i\pi x}{L}dx.  
\end{equation}
In the rest of this chapter we will focus on profiles $n(x)$ (and coefficients $\alpha_i$) having in mind that the analysis of profile $c(x)$ is done in the same way. 

\begin{figure}[h]
\begin{subfigure}{.32\textwidth}
\includegraphics[scale=0.22]{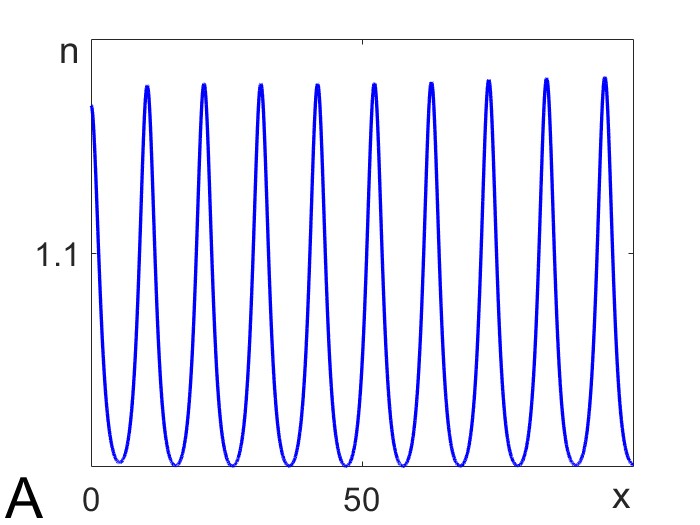}
\end{subfigure}
\begin{subfigure}{.32\textwidth}
\includegraphics[scale=0.22]{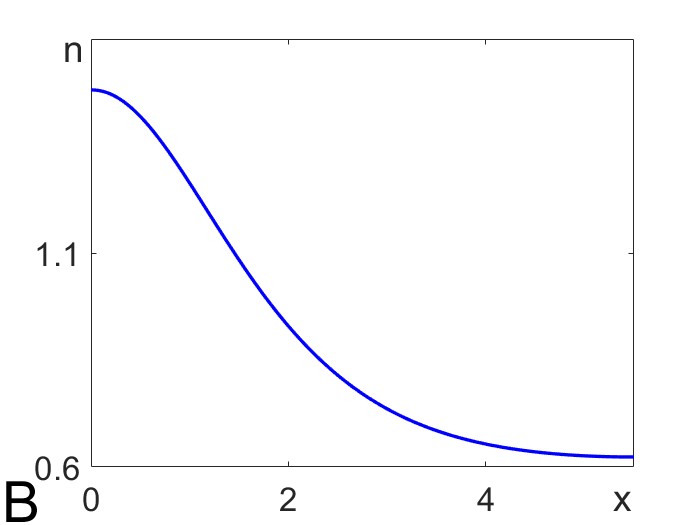}
\end{subfigure}
\begin{subfigure}{.32\textwidth}
\includegraphics[scale=0.22]{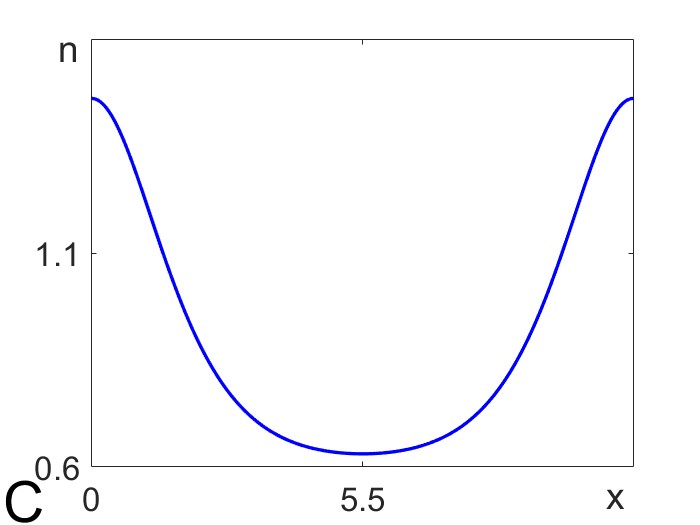}
\end{subfigure}
\caption{\em{\textbf{Profiles of cell density, $n(x)$, for Turing patterns in system \eqref{nondimensional} obtained numerically}. Domain size: $L=100$ in \textbf{A}, $L=5.5$ in \textbf{B}, and $L=11$ in \textbf{C}. Values of model parameters  $D=1$, $\chi_0=1.9$ and $r=0.1$.}}
\label{numerical_spikes}
\end{figure}

Formulas \eqref{alphas} can be used for spectral analysis of Turing patterns obtained numerically. A typical stationary profile for cell density, $n$, which is obtained in numerical simulations of the system \eqref{nondimensional} is shown in Fig.\ref{numerical_spikes}A. Spectral decomposition of this profile shows that only three modes have reasonably high coefficients, namely, $\alpha_0=0.8979$,  \mbox{ } $\alpha_{19}=0.3945$ and $\alpha_{38}=0.1293$, while all other coefficients are considerably smaller (lesser than 0.1). $\alpha_0$ gives the average level of $n$ for the entire pattern, $\alpha_{19}$ defines the amplitude of the mode whose characteristic length is 1/19th of the domain size, that is the one represented by 9.5 spikes. This is exact number of spikes we see in Fig.\ref{numerical_spikes}A. Finally, $\alpha_{38}$ corresponds to the second harmonic of the main harmonic given by $\alpha_{19}$. Thus,  the amplitude of the pattern shown in Fig.\ref{numerical_spikes}A can be estimated as $2\alpha_{19}$. Varying the size of the domain will result to a change in the number of observed spikes, however their amplitude and spatial periodicity wouldn't change. This is illustrated by patterns shown in panels B and C of Fig.\ref{numerical_spikes}, which are obtained in simulations of small domains where only half (B) and one (C) spike can fit in.

\begin{figure}[h]
\centering
\includegraphics[scale=.37]{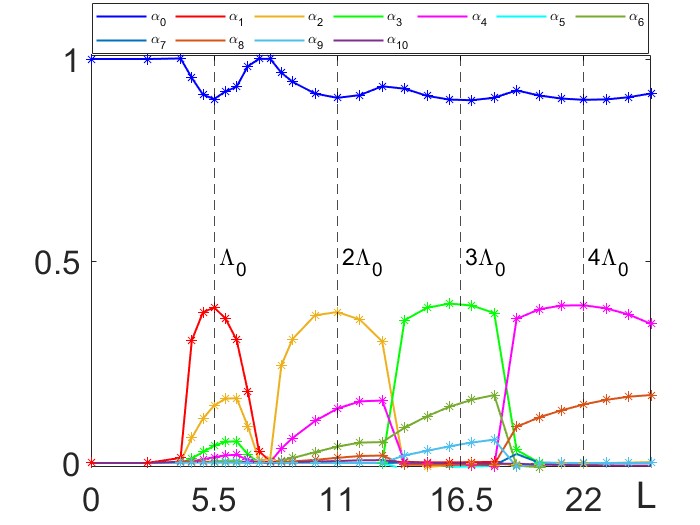}
\caption{\em{\textbf{Dependence of coefficients $\alpha_i$, $i=0, 1, ..., 10$ on domain size, $L$.} Plots are obtained by applying formula \eqref{alphas} to $n-$profiles from numerical simulations. Model parameters: $D=1$, $\chi_0=1.9$ and $r=0.1$. Vertical dashed lines indicate domain sizes given in whole number of characteristic length $\Lambda_0$.}}
\label{alphas_num}
\end{figure}

Patterns occurring in a small domain (like those shown in Fig.\ref{numerical_spikes}B and C) will be of particular interest for us, as corresponding Fourier series can be truncated at reasonably law value of $M$ (see eq. \eqref{sol}). Using Fourier decomposition of numerically obtained patterns we have investigated how the coefficients of Fourier modes depend on the size of domain. In Fig.\ref{alphas_num} we have shown plots of $\alpha_i$, for $i=0, 1, ...,10$, against the size of domain, $L$, which varied from $0$ to $25$. We can see that while the value of $\alpha_0$ does not change vary much (stays in the range between $0.9$ and $1$), and $\alpha_{10}$ is always negligibly small, other coefficients vary significantly. For $L<4$ $\alpha_0=1$ and all other coefficients are zero, which means that the system is in homogeneous state and no spikes have formed. For $4<L<7.5$ the first coefficient, $\alpha_1$, is larger than any subsequent coefficient ($\alpha_1>\alpha_i$ for any $i>1$), and this reflects the fact that only half of the spike can fit in. Within this range of domain sizes the value of $\alpha_1$ increases from zero and then decreases to zero reaching its maximum value $\alpha_1(=\alpha_{max})=0.385$ at $L=5.5$. This domain size, which corresponds to the maximal value of the coefficient $\alpha_1$ will be considered as a proper characteristic length (or half-wavelength) of Turing pattern and will be denoted as $\Lambda_0$, i.e. $\Lambda_0=5.5$ for the system \eqref{nondimensional} with model parameters in simulations used to produce Fig. \ref{alphas_num}. 

We see in Fig.\ref{alphas_num} that with the further increase of medium size ($L>7.5$) the coefficient $\alpha_2$ increases and then decreases, and then all the same with the coefficients $\alpha_3$ and $\alpha_4$, each reaching the maximum of $\alpha_{max}=0.385$ when the size of the domain $L=i*\Lambda_0$ for $\alpha_i$ (where $i=2$, $3$ and $4$),  which are indicated by vertical dashed lines. From these observations we conclude that Turing patterns observed in the system \eqref{nondimensional} with values of model parameters used to produce plots in Fig.\ref{alphas_num}  have the characteristic length $\Lambda_0=5.5$ (size of half-spike) and the amplitude $A=2\alpha_{max}=0.77$ (compare with the  amplitude of all three patterns shown in Fig. \ref{numerical_spikes}). 

In this section we have performed Fourier decomposition of patterns obtained in numerical simulations. Our next task is to find the Fourier coefficients describing Turing pattern analytically. This would let us to compare analytical results with numerical ones and besides predict the dependence of the wavelength, $\Lambda_0$ , and amplitude, $A$, of Turing patterns on model parameters in \eqref{nondimensional}.


\section{Nonlinear analysis of Turing patterns}

Stationary solutions of the system \eqref{nondimensional} represented by their Fourier series, satisfy the system \eqref{nondimensional} where partial derivatives over time are set to zero and variables $n$ and $c$ are replaced by their Fourier series \eqref{sol}. This gives the following equations for the coefficients of Fourier series, $\alpha_i$ and $\beta_i$: 
\begin{equation}
\begin{cases}
\displaystyle -D\sum_{i=0}^M (i\kappa)^2\alpha_i\cos(i\kappa x)+\chi_0\frac{\partial}{\partial x}\left(\sum_{i=0}^M \alpha_i\cos(i\kappa x)\sum_{i=0}^M i\kappa\beta_i\sin(i\kappa x)\right)+\\
\hfill r\sum_{i=0}^M \alpha_i\cos(i\kappa x)\left(1-\sum_{i=0}^M\alpha_i\cos(i\kappa x)\right)=0,\\[10pt]
\displaystyle -\sum_{i=0}^M(i\kappa)^2\beta_i\cos(i\kappa x)+\sum_{i=0}^M\alpha_i\cos(i\kappa x)-\sum_{i=0}^M \beta_i\cos(i\kappa x)=0,
\end{cases} \label{fourier_expansion}
\end{equation}
where $\kappa=\pi/L$. We will start analysis of this system by finding the homogeneous solutions, or the  solutions for which $\alpha_i=\beta_i=0$ when $i>0$. This means that in order to find the homogeneous solution we truncate the above system at $M=0$ and get: 
\begin{equation*}
\begin{cases}
\alpha_0-\alpha_0^2=0,\\[7pt]
\alpha_0-\beta_0=0.
\end{cases}
\end{equation*}
Roots of this system correspond to two homogeneous solutions: $\alpha_0=\beta_0=1$ (meaning that $u(x,t)=c(x,t)=1$, $\forall x,t$) and $\alpha_0=\beta_0=0$ (meaning that $u(x,t)=c(x,t)=0$, $\forall x,t$). 

To derive equations for $\alpha_i$ and $\beta_i$ when we truncate the system \eqref{fourier_expansion} at arbitrary $M$ we use trigonometric identities (to replace product of cosines by cosines of sums). By equating terms containing $\cos(ikx)$ for each $i$ we break the first equation in the system \eqref{fourier_expansion} into $M+1$ equations: 
\begin{equation}
\begin{cases}
\displaystyle \alpha_0-\alpha_0^2-\sum_{i=1}^M \frac{\alpha_i^2}{2}=0,\\[12pt]
\displaystyle -D\alpha_1k^2+\chi_0 k^2\left(\alpha_0\beta_1+\sum_{i \geq 1}^{M} \frac{i+1}{2}\alpha_i\beta_{i+1}-\sum_{i \geq 2}^M \frac{i-1}{2} \alpha_i\beta_{i-1}\right)+\\
\displaystyle \hfill r\left(\alpha_1-2\alpha_0\alpha_1-\sum_{i \geq 1}^M \alpha_i\alpha_{i+1}\right)=0,\\[12pt]
\displaystyle -4D\alpha_2k^2+\chi_0 k^2 \left( 4\alpha_0\beta_2+\alpha_1\beta_1+\sum_{i \geq 1}^M (i+2)\alpha_i\beta_{i+2} -\sum_{i \geq 3}^M (i-2)\alpha_i\beta_{i-2} \right)+\\
\displaystyle \hfill r\left( -\frac{\alpha_1^2}{2}+\alpha_2-2\alpha_0\alpha_2-\sum_{i \geq 1}^M \alpha_i\alpha_{i+2} \right)=0,\\[12pt]
\displaystyle -9D\alpha_3k^2 +\chi_0 k^2 \left(  9\alpha_0\beta_3 +3\alpha_1\beta_2+\frac{3}{2}\alpha_2\beta_1 +\sum_{i \geq 1}^M \frac{3(i+3)}{2}\alpha_i\beta_{i+3}-\sum_{i\geq 4}^M \frac{3(i-3)}{2}\alpha_i\beta_{i-3}  \right)+\\
\displaystyle \hfill r\left(\alpha_3-2\alpha_0\alpha_3-\alpha_1\alpha_2-\sum_{i \geq 1}^M \alpha_i\alpha_{i+3}    \right)=0,\\[12pt]  
....\\
\end{cases}\label{mode_expansion}
\end{equation}
Here the first equation is build as an identity for coefficients of $\cos(0)$, the second equation - for $\cos(kx)$, $i$-th equation - for $\cos(ikx)$. As the second equation in \eqref{fourier_expansion} is linear it breaks into $M+1$ identical equations: 
\begin{equation*}\label{cond1}
\displaystyle -\beta_i(ik)^2+\alpha_i-\beta_i=0.
\end{equation*}
This gives $\displaystyle \beta_i=\frac{\alpha_i}{1+(ik)^2}$ so that the coefficients $\beta_i$ can be taken out from the system \eqref{mode_expansion}, and we get a system of $M+1$ equations for $M+1$ unknowns. 

By truncating the system \eqref{mode_expansion} at $M=1$ and using the formula $\displaystyle \beta_1=\frac{\alpha_1}{1+k^2}$ we get the following system for finding coefficients $\alpha_0$ and $\alpha_1$: 
\begin{equation}
\begin{cases}
\displaystyle \alpha_0-\alpha_0^2-\frac{\alpha_1^2}{2}=0,\\[7pt]
\displaystyle -Dk^2\alpha_1+k^2\chi_0\alpha_0\frac{\alpha_1}{1+k^2}+r(-2\alpha_0\alpha_1+\alpha_1)=0.
\end{cases} \label{m1_system}
\end{equation}
This system has four solutions of which two, namely $(0, 0)$ and $(1,0)$, correspond to  homogeneous solutions of the system \eqref{nondimensional} and two others,
\begin{equation}\label{a0}
\displaystyle \alpha_0= \frac{(Dk^2-r)(1+k^2)}{\chi_0 k^2-2r(1+k^2)} \mbox{ } \mbox{and} \mbox{ } \displaystyle \alpha_1=\pm  \sqrt{2\alpha_0(1-\alpha_0)}
\end{equation}
- to spatial periodic patterns forming in \eqref{nondimensional}. The roots associated with homogeneous solutions were found earlier after truncating the system \eqref{mode_expansion} at $M=0$ while the solutions \eqref{a0} correspond to the first mode of the spectral decomposition of Turing pattern. Thus, equations \eqref{a0} give the dependence of $\alpha_0$ and $\alpha_1$ on model parameters and can be compared with the results obtained numerically and shown in Figure \ref{alphas_num}.

Truncation of the system \eqref{mode_expansion} at $M=2$ gives three equations for three unknowns, $\alpha_0$, $\alpha_1$ and $\alpha_2$:
\begin{equation*}
\begin{cases}
\displaystyle \alpha_0-\alpha_0^2-\frac{\alpha_1^2}{2}-\frac{\alpha_2^2}{2}=0,\\[7pt]
\displaystyle -D\alpha_1k^2+\chi_0 k^2\left(\alpha_0\beta_1+\alpha_1\beta_2-\frac{1}{2}\alpha_2\beta_1\right)+ r(\alpha_1-2\alpha_0\alpha_1-\alpha_1\alpha_2)=0,\\[7pt]
\displaystyle -4D\alpha_2k^2+\chi_0 k^2(4\alpha_0\beta_2+\alpha_1\beta_1)+r\left(\alpha_2-\frac{\alpha_1^2}{2}-2\alpha_0\alpha2\right)=0.
\end{cases}
\end{equation*}
Similarly, systems of equations for 4 and 5 unknowns, when the system \eqref{mode_expansion} is truncated at $M=3$ or $M=4$, can be derived. Although roots of these system can't be expressed explicitly, their values (for a given set of model parameters) can be found using Matlab code. Using such code we have find solutions of the system \eqref{mode_expansion} truncated up to $M=4$. These solutions for the model \eqref{nondimensional} with values of model parameters $D=1$, $\chi_0=1.9$ and $r=0.1$ (see Figs \ref{numerical_spikes} and \ref{alphas_num}) and domain size $L=5.5$ are shown in Table \ref{coefficients_table}. For comparision we have also shown in this table the coefficients $\alpha_i$ obtained by spectral decomposition of the numerically simulated pattern (shown in Fig.\ref{numerical_spikes}B). As we can see from this table, the obtained solutions are consistent across different truncations and fairly comparable with the results of spectral decomposition of simulated pattern. Also, the coefficients $\alpha_3$ and $\alpha_4$ are considerably small and this indicates that the truncation of system \eqref{mode_expansion} at $M=4$ should give a fairly accurate solution for the model \eqref{nondimensional} when the domain size is reasonably small (in the presented case $L=\Lambda_0$). 

\begin{table}
	\begin{center}
		\begin{tabular}{| l | l |  l |  l | l | l | } 
			\hline
			$ $  & $\alpha_0$ & $\alpha_1$ & $\alpha_2$ & $\alpha_3$ & $\alpha_4$ \\[7pt] \hline
			Simulation & 0.9009 & 0.3849 &  0.1423 & 0.0428 & 0.0143 \\[7pt] \hline
			$M=1$ & 0.8462 & 0.5103 &0 & 0 & 0 \\[7pt] \hline
			$M=2$ & 0.8868 & 0.4153 & 0.1680 & 0 & 0 \\[7pt] \hline
			$M=3$ & 0.8938 & 0.4035 & 0.1563 & 0.0508 & 0\\[10pt] \hline
			$M=4$ & 0.8941 & 0.4033 & 0.1550 &0.0496 & 0.0155 \\[10pt] \hline
		\end{tabular}
	\end{center}
	\caption{\em{\textbf{Coefficients of Fourier series, $\alpha_i$, for the first five modes of Turing pattern represented by half-spike.} Coefficients in the first row obtained using eq. \eqref{alphas} for the profile shown in Fig. \ref{numerical_spikes}. Values in the following four rows obtained by truncation of the system \eqref{mode_expansion} at $M=1$, $2$, $3$ and $4$. Values of the model parameters are the same as in Fig. \ref{numerical_spikes}B.}}
	\label{coefficients_table}
\end{table}

\begin{figure}[h]
	\centering
	\begin{subfigure}[b]{0.48\textwidth}
\includegraphics[scale=.32]{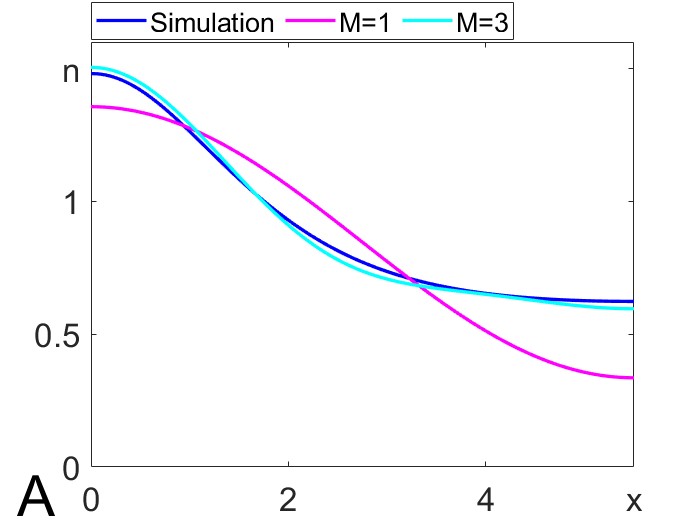}
	\end{subfigure}
	\begin{subfigure}[b]{0.48\textwidth}
\includegraphics[scale=.32]{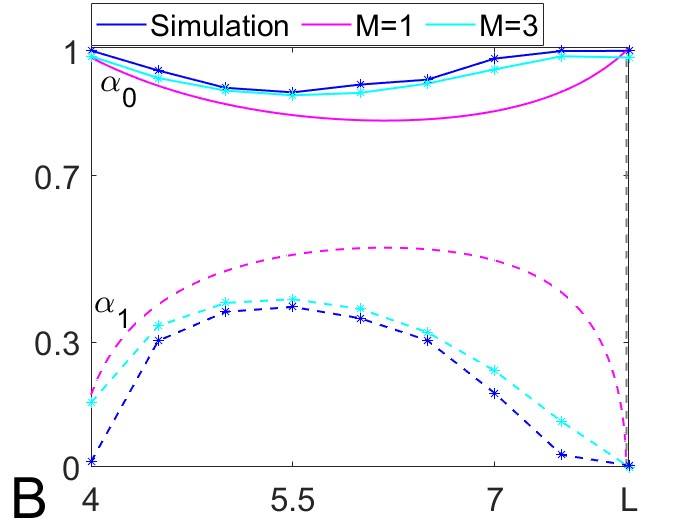}
	\end{subfigure}
\caption{\em{\textbf{Comparison of analytical and numerical solutions.} \textbf{A}: Profiles of solutions to \eqref{mode_expansion} truncated at $M=1$ (red) and $M=3$ (cyan) as compared to numerically simulated profile (blue, same as in Fig. \ref{numerical_spikes}B). \textbf{B}: Dependence of the coefficients $\alpha_0$ (solid lines) and $\alpha_1$ (dashed lines) on the domain size $L$ found by spectral decomposition of profiles obtained in simulations (blue) and from \eqref{mode_expansion} truncated at $M=1$ (red) and $M=3$ (cyan). Values of model parameters are the same as in Fig.\ref{numerical_spikes}.}}
	\label{profile_comp1}
\end{figure}

To check the accuracy of truncated solutions presented in Table \ref{coefficients_table} we have reproduced corresponding profiles (for $M=1$, $2$ and $3$) and compared them with the profile obtained in numerical simulations. All these profiles are shown in Fig. \ref{profile_comp1}A. Visually it is evident that the solution of \eqref{mode_expansion} truncated at $M=2$ fits better  the profiles obtained in numerical simulation than the one truncated at $M=1$. Furthermore, the profile corresponding to the solution truncated at $M=3$ is almost indistinguishable from the simulated profile. As a numerical measure of discrepancy between profiles, $n_i(x)$ and $n_0(x)$, one can use the integral:  
\begin{equation*}
I=\int_0^L [n_i(x)-n_0(x)]^2 dx,
\end{equation*} \label{discrep}
where $n_i(x)$ correspond to one of the solutions of \eqref{mode_expansion} (truncated at $M=i$, $i>0$), while $n_0(x)$ - to the profile obtained in numerical simulations. The value of this integral for the solutions truncated at $M=1$ is $I=0.1353$, at $M=2$ - $I=0.0114$ and at $M=3$ - $I=0.0026$. Thus, the discrepancy between the solution, truncated at $M=3$, and the one obtained numerically is lesser than 1\%, and besides the numerical discrepancy for $M=4$ is $I=0.0020$ and fairly the same as for $M=3$. 

The amplitude of the profile for variable $n$ in a small domain (when $L \approx \Lambda_0$) is predominantly defined by the coefficient $\alpha_1$, while its average value by $\alpha_0$. Figure \ref{profile_comp1}B shows the dependence of these two coefficients on the size of the domain, $L$, obtained from the spectral decomposition of the simulated profile and according to the equation \eqref{fourier_expansion} truncated at $M=1$ and $M=3$. We see that while all three plots are qualitatively similar, the plots for simulated profile and for the solution of eq. \eqref{fourier_expansion} truncated at $M=3$ also show a good quantitative agreement. 

In this section we have focused on finding the analytical solutions of the system \eqref{nondimensional} represented by truncated Fourier series. Using Matlab code we were able to find solutions truncated up to $M=4$, and have shown that these solutions (at $M \ge 3$) are fairly accurate if the modelled domain is small. Throughout this section we have fixed the model parameters: $D=1$, $\chi_0=1.9$ and $r=0.1$. In the next section we will check how the characteristics of Turing patterns obtained numerically and given by the analytical solutions of \eqref{fourier_expansion} depend on the model parameters.


\section{Impact of model parameters on characteristics of Turing pattern}
\begin{figure}[h]
	\centering
	\begin{subfigure}[b]{0.33\textwidth}
		\centering
		\includegraphics[width=\textwidth]{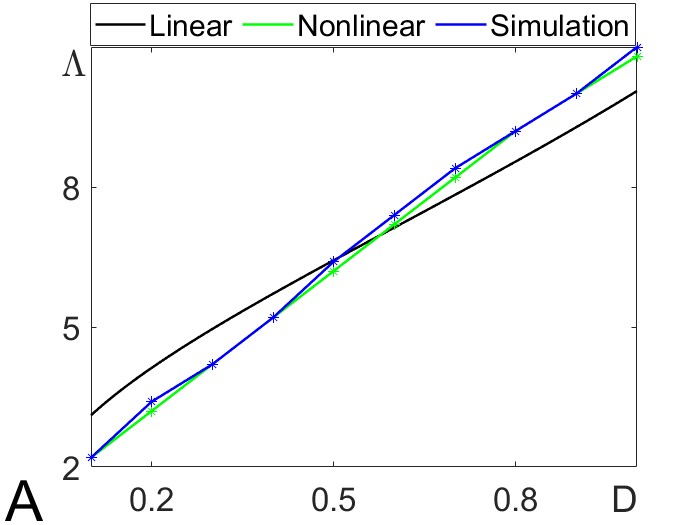}
	\end{subfigure}
	\begin{subfigure}[b]{0.325\textwidth}
		\centering
		\includegraphics[width=\textwidth]{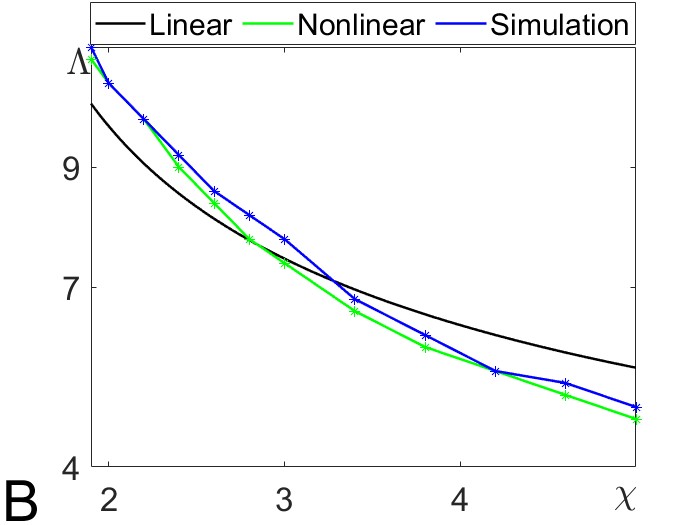}
	\end{subfigure}
\begin{subfigure}[b]{0.33\textwidth}
		\centering
		\includegraphics[width=\textwidth]{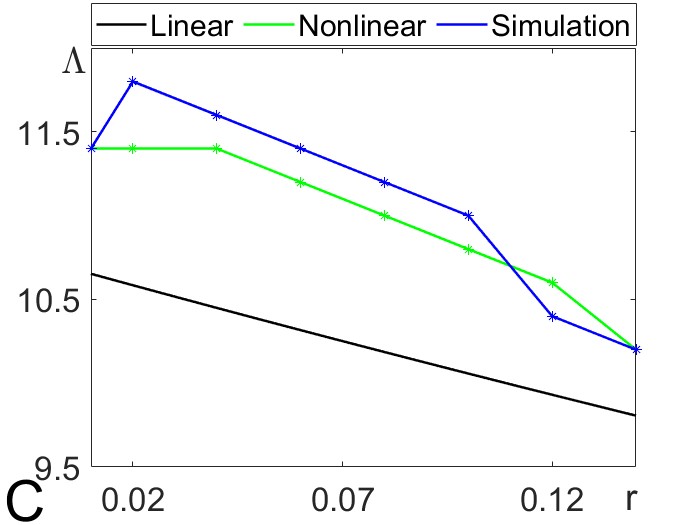}
	\end{subfigure}
	\begin{subfigure}[b]{0.33\textwidth}
	\centering
	\includegraphics[width=\textwidth]{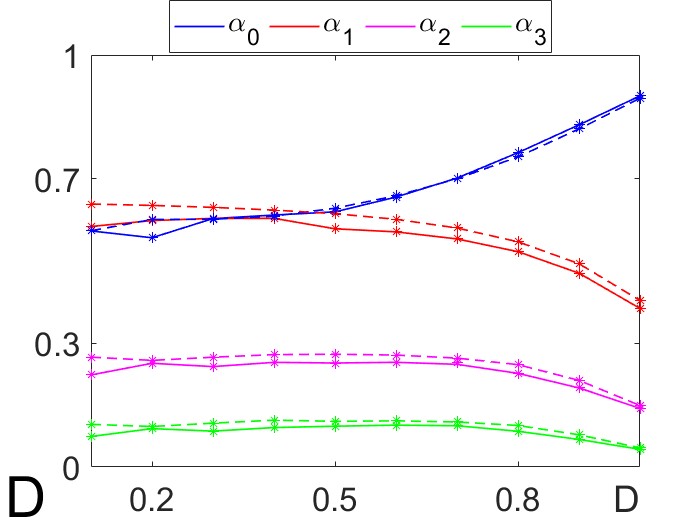}
\end{subfigure}
\begin{subfigure}[b]{0.325\textwidth}
	\centering
	\includegraphics[width=\textwidth]{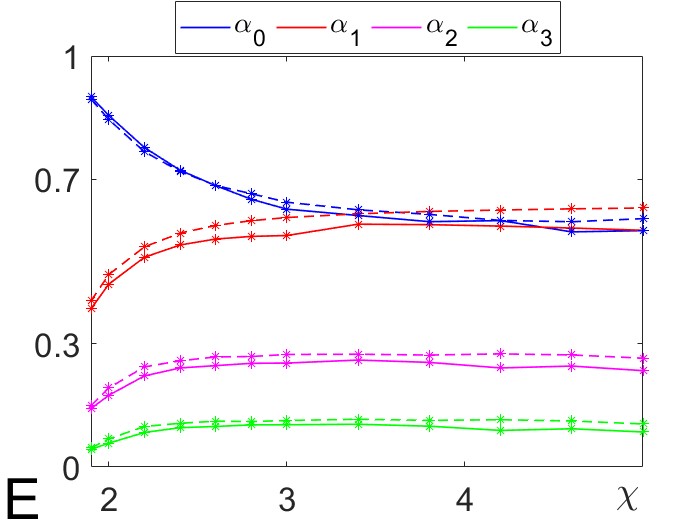}
\end{subfigure}
\begin{subfigure}[b]{0.33\textwidth}
	\centering
	\includegraphics[width=\textwidth]{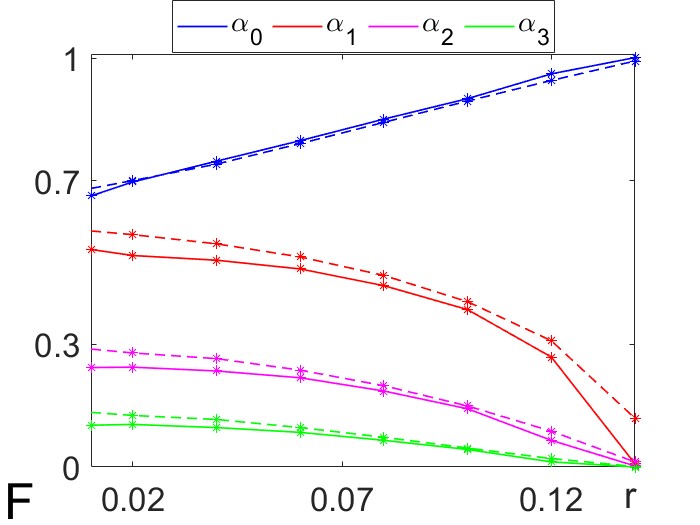}
\end{subfigure}
	\begin{subfigure}[b]{0.33\textwidth}
	\centering
	\includegraphics[width=\textwidth]{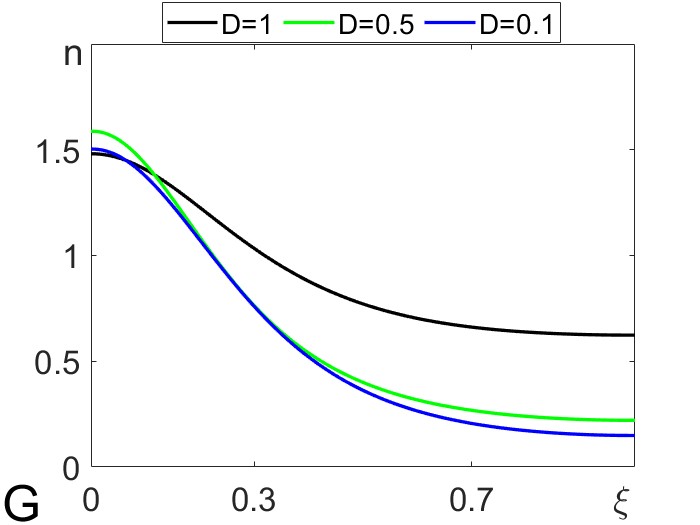}
\end{subfigure}
\begin{subfigure}[b]{0.325\textwidth}
	\centering
	\includegraphics[width=\textwidth]{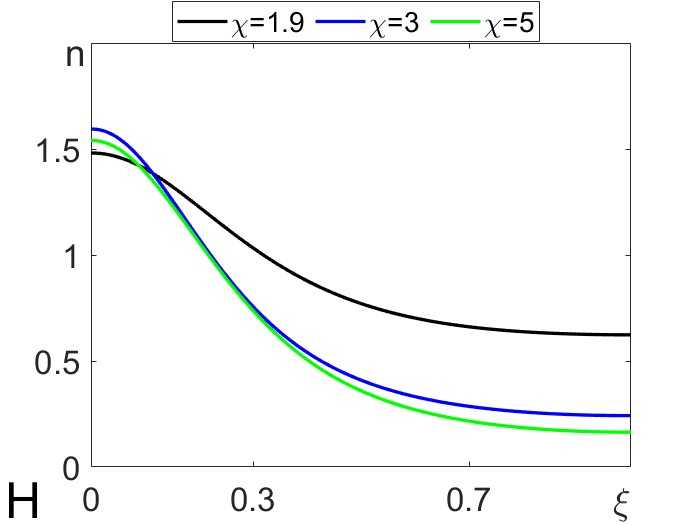}
\end{subfigure}
\begin{subfigure}[b]{0.33\textwidth}
	\centering
	\includegraphics[width=\textwidth]{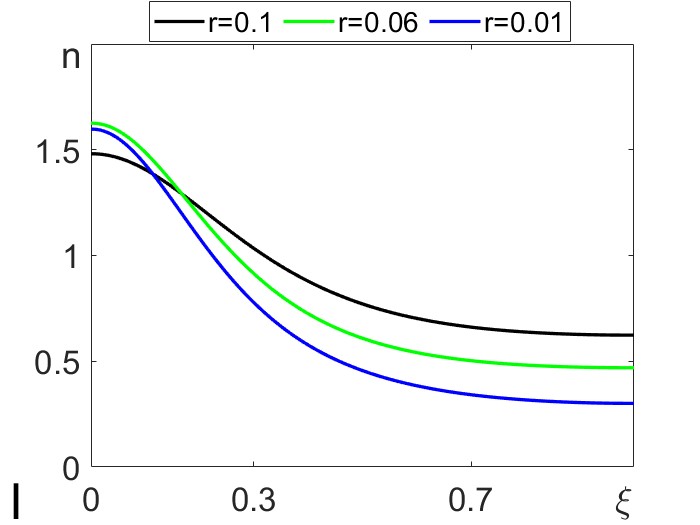}
\end{subfigure}
\caption{\em{\textbf{Dependence of the characteristics of Turing pattern on model parameters.} \textbf{A-C} Dependence of the wavelength on the diffusion coefficient, $D$, (panel \textbf{A}), chemotactic coefficient, $\chi_0$, (panel \textbf{B}), and reproduction rate, $r$, (panel \textbf{C}): black lines - from linear analysis (same as solid line in Fig. \ref{linear_approx}), green lines - from Fourier series truncated at $M=3$, and blue lines - from simulations. \textbf{D-F} Coefficients of Fourier series versus the diffusion coefficient, $D$, (panel \textbf{D}), chemotactic coefficient, $\chi_0$, (panel \textbf{E}), and reproduction rate, $r$, (panel \textbf{F}) found after spectral decomposition of numerically simulated profile (solid lines) and found analytically from eq. \eqref{mode_expansion} with $M=3$ (dashed lines). \textbf{G-I} Profiles of cell density obtained numerically for different values of the diffusion coefficient, $D$, (panel \textbf{G}), chemotactic coefficient, $\chi_0$, (panel \textbf{H}), and reproduction rate, $r$, (panel \textbf{I}). Each profile obtained in the domain of size of characteristic length, $L=\Lambda_0$ (for given set of parameter values) and presented versus scaled spatial variable $\xi=x/\Lambda_0$. }}
\label{par}
\end{figure}

Based on the linear analysis of pattern formation presented in Section 3 we were able to predict dependence of the pattern's wavelength on model parameters (see eqs \eqref{lamb1}, \eqref{lamb2} and Fig. \ref{linear_approx}). However, linear analysis didn't let us to make any predictions concerning the amplitude of Turing pattern. In this section we will make such predictions using nonlinear analysis, presented in the previous section, and compare them with results of numerical simulations. The outcome of this study is presented in Fig. \ref{par}. Since $D$ represents the ratio of the rate of random motion of cells to the diffusion coefficient of chemotactic agent it is essentially smaller than one and therefore $D<1$ in our plots. Also, ranges of values for $\chi_0$ and $r$ were chosen in a way that the condition for Turing instability given by the equation \eqref{m3_rt} (and illustrated in Fig. \ref{ratio_implicit}) is satisfied, that is $\chi_0>1.9$ or $r<0.14$ for the default set of parameters.

The dependence of pattern's wavelength on model parameters as predicted on the basis of linear analysis (eq. \eqref{lamb2}), nonlinear analysis (eq. \eqref{mode_expansion}) and found numerically is shown on panels A-C of Fig. \ref{par}. The wavelength in numerical simulations was found as a length of the domain for which Fourier decomposition of obtained profile gave maximal value of $\alpha_2$ (or $2\Lambda_0$ according to Fig. \ref{alphas_num}). Similarly, using equation \eqref{mode_expansion} truncated at $M=4$ we found the value of $k$ for which the value of $\alpha_2$ is maximal, and the wavelength was found as a length of domain corresponding to this value of $k$, namely, $\Lambda=\pi/k$. As we can see from Fig. \ref{par}A-C the dependences found this way indicate that the wavelength is increasing (linearly) with the diffusion coefficient, $D$, and decreasing with chemotactic coefficient, $\chi_0$, (can be approximated by the hyperbola) and reproduction rate, $r$, (appears as a linear) as it was predicted by linear analysis in Section 3. Concerning the dependence on the reproduction rate, $r$, the prediction made by the implicit formula \eqref{lamb2} rather than by the explicit \eqref{lamb1} (see Fig. \ref{linear_approx}) was confirmed both by simulations and analytically. It is evident that the predictions made using nonlinear analysis (on the basis of equation \eqref{mode_expansion}) are far more accurate as compared to those made using linear analysis (on the basis of equation \eqref{lamb2}) in reproducing numerical results.   

The dependence of the amplitude of Turing pattern on model parameters is addressed on panels D-I of Fig. \ref{par}. The amplitude is predominantly defined by the maximal value of $\alpha_i$ ($i>0$) when the size of the domain $L$ is varied (see Fig. \ref{alphas_num}). As it was found in Section 5, for the default values of model parameters, the maximal value of $\alpha_i$ is $\alpha_{max}=0.385$ for $i=1$, $2$, $3$, $4$ and evidently for all other modes. So, the main question we address now is how $\alpha_{max}$ depends on the model parameters. Panels D-F show how $\alpha_{max}$ which is given here by the value of $\alpha_1$ depends on model parameters. In addition we have presented here the values of $\alpha_0$, $\alpha_2$ and $\alpha_3$ which contribute to the shape of the $n-$profile. We do not consider contributions of $\alpha_2$ and $\alpha_3$ to the amplitude of Turing pattern for the following reasons: (1) due to symmetry in Fourier expansion, even modes do not contribute to the amplitude at all and therefore $\alpha_2$ has no impact to the amplitude of pattern; (2) according to Table \ref{coefficients_table} $\alpha_3$ by order of magnitude smaller than $\alpha_1$ which means that it contributes only about 10\% to the amplitude of patterns.  Besides, the plots in panels D-F indicate that this contribution is roughly the same for any set of model parameters and thus we can consider $\alpha_1$ as an indicator for the behaviour of the pattern's amplitude. 

The main conclusion drawn from the plots presented on panels D-F of Fig. \ref{par} is that the amplitude of the Turing pattern (given by the value of $\alpha_1$) saturates at low values of the diffusion, $D$, and reproduction rate, $r$, and at high values of the chemotactic sensitivity, $\chi_0$, while it increases with the increase of $D$ and $r$ and decreases with the increase of $\chi_0$. Another notable conclusion, one can draw from the plots presented on these panels, is that the average value of cell density, $n$, given by the value of $\alpha_0$, is increasing with the diffusion coefficient $D$ and reproduction rate, $r$, while decreasing with the chemotactic sensitivity, $\chi_0$. It saturates at low values of $D$ and high values of $\chi_0$ while depends linearly on $r$.  

Panels G-I of Fig. \ref{par} show how the shape of $n-$profile depends on model parameters. Plots on panel G confirm that the amplitude of the profile decreases and the average value of the density of cells, $n$, increases with the increase of diffusion coefficient, $D$. Plots on panel H confirm that the amplitude of the profile increases and the average value of the density of cells, $n$, decreases with the increase of chemotactic sensitivity, $\chi_0$. And finally, plots on panel I confirm that the amplitude of the profile decreases and the average value of the density of cells, $n$, increases with the increase of reproduction rate, $r$. Dependence of the wavelength of patterns on model parameters can't be seen on panels G-I as the profiles there, although obtained for media of different sizes, are spatially scaled and presented as $n(x)=n(\xi)$ in scaled variable $\xi=x/L$.       


\section{Discussion}

In this paper, we have analysed periodic stationary patterns forming in a model of growing population of motile bacteria. The growth of modelled population is following the logistic law while migration of cells is due to chemoattraction by the chemical produced by the cells themselves. In terms of mathematics, we have studied stationary solutions of nonlinear two variable reaction-diffusion-advection system. It is known that nonlinear system without advective term should have at least cubic terms to produce periodic patterns of finite amplitude \cite{Vasiev2013, Vasiev2016a}. It is also known that for formation of periodic patterns the advection term should represent chemoattraction \cite{Ks}, while chemorepulsion only causes formation of moving bands of cells \cite{Keller1, Vasiev2011, Vasiev2016}. 

In our analysis of the models of type \eqref{nondimensional} we have identified conditions for formation of Turing patterns (or for Turing instability) for a number of commonly used models \cite{Painter}. This conditions where presented as a certain expression, which we denoted as $R_T$, given in terms of model parameters and such that $R_T>1$ for Turing pattern to emerge. The expressions defining  $R_T$ for a number of models are presented in Table \ref{tab_models}. It appears that the domains in parameter space defined by the condition $R_T>1$ for different models are very similar, which is illustrated in Fig. \ref{ratio_implicit} for models $M3$ and $M7$. The main conclusion of this analysis is that for Turing patterns to form, the chemoattraction should be strong enough, while the diffusion of chemical and the reproduction of cells are slower than certain threshold rate defined by the condition $R_T>1$. 

We have also presented our analysis of spatial characteristics of Turing patterns. Commonly, the wavelength of Turing pattern is estimated on the basis of linear analysis using formula \eqref{lamb1} \cite{Murray2}. We have suggested a correction to this estimate, which is also based on linear analysis and given by the implicit formula \eqref{lamb2}. Comparison with the wavelengths of simulated patterns has indicated that the formula \eqref{lamb2} gives better approximation to the wavelength of Turing pattern not only quantitatively but also qualitatively when it comes to checking the dependence of the wavelength on model parameters (see Fig. \ref{linear_approx} and Fig. \ref{par}, panel C). The wavelength of Turing pattern was also estimated on the basis of nonlinear analysis presented in Section 6. The first order approximation can be made analytically by finding the value of $k$ corresponding to the maximal value of $\alpha_1$ in equations \eqref{a0}. The link between this value of $k$ with the wavelength is given by formula $k=\pi/\Lambda_0$ where $\Lambda_0$ is a characteristic length, which is twice smaller than the wavelength. Dependence of $\alpha_1$ on the domain size, $L$ found on the basis of this approximation is shown by dashed-red plot in Fig. \ref{profile_comp1}, panel B. It is evident from this figure that the length of the medium for which $\alpha_1$ reaches its maximal value (and which gives the value of $\Lambda_0$) is much higher that the one found numerically ($\Lambda_0=5.5$, see Fig. \ref{alphas_num}). However when we find the characteristic length using the model \eqref{mode_expansion} truncated at $M=3$ we find that it is very close to the one obtained numerically (compare locations of maxima for dashed blue and cyan in Fig. \ref{profile_comp1}, panel B).  

One of the most important problems addressed in this work is about analytical estimation of the amplitude of Turing pattern. This problem gets recently attention of many researchers. Using multi-scale analysis by means of amplitude equation for analysis of Turing patterns turned to be technically sophisticated and factually not convincing \cite{Dutt,Wang}, although further development of this method may turn to be more useful. The linear approach for finding the amplitude of Turing pattern suggested in \cite{Buceta} looks somehow oversimplified. In terms of the system \eqref{nondimensional} it brings to the equation for $\alpha_1$ given by the second equation in system \eqref{m1_system} where the value of $\alpha_0$ is taken from the homogeneous solution, $\alpha_0=1$. We have found that the value of $\alpha_1$ found this way fits numerical results worth than the one found as a solution of the system \eqref{m1_system}. Furthermore, the amplitude found from system \eqref{mode_expansion} truncated at $M=3$ is far better even when compared with the one found using system \eqref{m1_system} and the Turing pattern reproduced using this solution is almost indistinguishable from the one obtained in numerical simulations (see Fig. \ref{profile_comp1}, panel A). The only concern about finding solutions of algebraic system \eqref{mode_expansion} is that this solution can be expressed explicitly only when this system truncated at $M=1$ (see \eqref{a0}), while when $M>1$ solutions can be found only numerically using Matlab (or Maple) code. This is especially true about the dependence of wavelength and amplitude of Turing patterns found for the system \eqref{mode_expansion} and presented in Fig. \ref{par}. All corresponding plots are drawn on the basis of sets of solutions of system \eqref{mode_expansion} truncated at $M=3$, which where found numerically.

In this paper we have presented the analysis of spatial characteristics and amplitude of Turing patterns forming in the model given by system \eqref{nondimensional}. As this model is a representative for a wide class of models (see Table \ref{tab_models}), presented results, at least, qualitatively are applicable to many other models. The methods used for the presented analysis, particularly, representation of Turing patterns by means of Fourier series, appeared to be very helpful. It can be used for the analysis of other classes of models where formation of Turing patterns can be observed. Particularly, it can be extended for the analysis of patterns forming in a system of two competing species of which one is producing chemical agent and the other responds to it chemotactically. We are planning to present such study in the follow up publication.  

\section*{Acknowledgment}
This work has been supported by the EPSRC (Engineering and Physical Sciences  Research Council) and NBIC (National Biofilms Innovation Centre) scholarships to Valentina Bucur. Authors are grateful to Prof Rachel Bearon for helpful discussions.


\end{document}